\documentclass{article}

\usepackage{url}

\usepackage[T1]{fontenc}     % Output font encoding for international characters
\usepackage[utf8]{inputenc}
\usepackage{csquotes}
\usepackage[a4paper]{geometry}
\usepackage{amsmath}
\usepackage{mathtools}
\usepackage{graphicx}
\usepackage[english]{babel}
\usepackage{hyperref}

\usepackage{tikz}
\usetikzlibrary{math}
\usetikzlibrary{decorations.pathreplacing}
\usepackage{placeins}
\usepackage{amssymb,amsfonts,amsthm}
\usepackage{dsfont}
\usepackage{xcolor}
\usepackage{todonotes}
\usepackage{appendix}
\usepackage{caption}
\usepackage{subcaption}

\usepackage[backend=biber, style=numeric]{biblatex}
\newtheorem{proposition}{Proposition}

\newcommand{\N}{\mathbb N}

\newcommand{\R}{\mathbb R}

\newcommand{\cH}{\ensuremath{\mathcal{H}}}

\newcommand{\bbc}{\boldsymbol{c}}

\newcommand{\bbm}{\boldsymbol{m}}

\newcommand{\bbF}{\boldsymbol{F}}

\newcommand{\bbJ}{\boldsymbol{J}}
\newcommand{\bbA}{\boldsymbol{A}}
\newcommand{\bbM}{\boldsymbol{M}}

\newcommand{\bbmu}{\boldsymbol{\mu}}

\newcommand{\ovc}{\bar{c}}
\newcommand{\ovX}{\bar{X}}

\renewcommand{\l}{\left}
\renewcommand{\r}{\right}

\mathtoolsset{showonlyrefs}

\begin{document}

\title{Structure Preserving Finite Volume Approximation of Cross-Diffusion Systems Coupled by a Free Interface}

\author{Clément Cancès, Jean Cauvin-Vila, Claire Chainais-Hillairet, Virginie Ehrlacher}

\date{}

\maketitle              % typeset the title of the contribution

\begin{abstract}
We propose a two-point flux approximation finite-volume scheme for the approximation of two cross-diffusion systems coupled by a free interface to account for one-dimensional vapor deposition. The moving interface is addressed with a cut-cell approach, where the mesh is locally deformed around the interface. The scheme preserves the structure of the continuous system, namely: mass conservation, nonnegativity, volume-filling constraints and decay of the free energy. Numerical results illustrate the properties of the scheme. 

\medskip

\noindent \textit{Keywords:} cross-diffusion system, cut-cell method, finite volume scheme, free energy dissipation, moving interface, vapor deposition
\end{abstract}

\section{A Free Interface Cross-Diffusion Model}
We address a toy model to describe a physical vapor deposition process used for the fabrication of thin film layers \cite{ccce2023}. We consider the evolving domain 
\[\Omega(t) = (0,X(t)) \cup (X(t),1), ~ t>0, \]
where $\R_+ \ni t \to X(t) \in [0,1] $ is the free interface between the solid (left) and the gas (right). Traces and jumps at the interface are respectively denoted by $f^s,f^{g}, [[f]] = f^g-f^s $. We consider $n$ different chemical species represented by their densities of molar concentration $\bbc = (c_1,\dots,c_n)^T$. The local conservation of matter reads:
\begin{subequations}
  \label{eq:model}
\begin{equation}
\partial_t \bbc + \partial_x \bbJ = 0, ~ t > 0, ~ x \in \Omega(t),
\label{eq:conservation}
\end{equation}
for some molar fluxes $\bbJ:= (J_1, \dots, J_n)^T$. Cross-diffusion phenomena are modelled differently in each phase. In the solid phase, the fluxes are given by 
\begin{equation}
  J_i = -\sum_{j=1}^n \kappa_{ij}^s \left(c_j \partial_x c_i - c_i \partial_x c_j \right), ~ \text{in} ~ (0,X), \; i \in \{1,\dots,n\},
  %\label{eq:size-exclusion-flux}
\end{equation}
with cross-diffusion coefficients $\kappa_{ij}^s=\kappa_{ji}^s > 0$, which rewrite more compactly
\begin{equation}
  \bbJ = - \bbA_s(\bbc) \partial_x \bbc, ~ \text{in} ~ (0,X),
  \label{eq:size-exclusion-flux}
\end{equation}
with a linear diffusion matrix $\bbA_s(\bbc)$ (see \cite{bakhta2018}). In the gaseous phase, the fluxes are defined implicitly via the Maxwell-Stefan linear system (see \cite{bothe2011}) 
\begin{equation}
    \bbA_g(\bbc) \bbJ = - \partial_x \bbc, ~ \text{and} ~ \sum_{i=1}^n J_i = 0, ~ \text{in} ~ (X,1),
    \label{eq:Maxwell-Stefan-flux}
\end{equation}
where $\bbA_g(\bbc)$ is identical to $\bbA_s(\bbc)$, except for possibly different cross-diffusion coefficients $\kappa_{ij}^g = \kappa_{ji}^g > 0$. The system is completed with an initial condition $(\bbc^0,X^0)$, no-flux conditions on the fixed boundary and the following conditions across the moving interface:
\begin{equation}
  \label{eq:interface-conservation}
  \bbJ^s(t) - X'(t) \bbc^s(t) = \mathds{1}_{\{X(t) \in (0,1)\}} \bbF(t) = \bbJ^g(t) - X'(t) \bbc^g(t), ~ t>0,
  \end{equation}
where $\bbF$ accounts for reaction mechanisms \cite{glitzky2013,mielke2014} and is defined, for some constant reference chemical potentials $\mu_i^{*,s}, \mu_i^{*,g} \in \R$, by the Butler-Volmer formulas: for $i \in \{1,\dots,n\}$,
\begin{equation}
  \label{eq:BV}
  \begin{aligned}
  F_i &=  c_i^s \exp \left(\frac{\mu_i^{*,g}-\mu_i^{*,s}}{2}\right) - c_i^g \exp \left(\frac{\mu_i^{*,s}-\mu_i^{*,g}}{2}\right),  \\
  &= 2 \sqrt{c_i^s c_i^g} \sinh\left(-\frac{1}{2} [[\log(c_i)-\mu_i^*]]\right).
  \end{aligned}
  \end{equation}
Finally, the interface evolves according to
\begin{equation}
\label{eq:interface}
X'(t) =  - \mathds{1}_{\{X(t) \in (0,1)\}} \sum_{i=1}^n F_i. 
\end{equation}
Note that, in the limit cases $X(t)=0 \text{ or } X(t)=1$, \eqref{eq:interface-conservation}-\eqref{eq:interface} imply that we recover a single phase problem with zero-flux boundary conditions.
\end{subequations}
The system enjoys several important properties we aim at preserving at the discrete level: First, mass conservation follows from the local conservation \eqref{eq:conservation}, no-flux conditions on the fixed boundary and the conservative condition \eqref{eq:interface-conservation}. Second, the system preserves the nonnegativity of the concentrations and the volume-filling constraints $\sum_{i=1}^n c_i =1$ (satisfied by the initial condition), and we refer to such a solution as \emph{admissible}. Finally, the functional
\begin{equation}
  \label{eq:free-energy}
  \cH(\bbc,X) = \int_0^X h_s(\bbc) + \int_X^1 h_g(\bbc),
\end{equation}
with density $h_{\alpha}(\bbc) = \sum_{i=1}^n c_i (\log(c_i)-\mu_i^{*,\alpha}) -c_i +1$, for $\alpha \in \{s,g\}$, can be shown to formally satisfy, for some positive semi-definite mobility matrices $\bbM_s,\bbM_g$, the free energy dissipation relation \cite{cances2020,cances2020d} 
\begin{equation}
  \label{eq:dissip}
  \begin{aligned}
&\frac{d}{dt} \cH(\bbc(t),X(t)) =  -\int_0^{X(t)} \partial_x \log(\bbc)^T \bbM_s(\bbc)\partial_x \log(\bbc) \\
&- \int_{X(t)}^1 \partial_x \log(\bbc)^T \bbM_g(\bbc)  \partial_x \log(\bbc)+ \bbF(t)^T [[\log(\bbc) - \bbmu^*]] \leq 0.
\end{aligned}
\end{equation}
One deduces from the dissipation inequality that stationary solutions $(\boldsymbol{\ovc},\ovX)$ must be constant in (each connected part of) $\bar{\Omega} = (0,\ovX) \cup (\ovX,1)$ and moreover, if $\ovX \in (0,1)$, $F_i(\ovc_i^s,\ovc_i^g)=0$ should hold for any $i$. We characterize in \cite{ccce2023} the stationary states of \eqref{eq:model}, as partially stated in Proposition~\ref{prop:stationary}. 
\begin{proposition}[Stationary states]
  \label{prop:stationary}
  Let $\bbm^0 = \bbm^{0,s} + \bbm^{0,g} > 0$ be the initial amount of matter in the system. The one-phase solutions $(\bbm^0,0,1)$ and $(0,\bbm^0,0)$ are stationary. Define the coefficients $\beta_i = \exp\l([[\mu_i^*]] \r)$. There exists a stationary solution where the two phases coexist (\emph{i.e.} such that $\ovX \in (0,1)$) if and only if  
  \begin{equation}
      \label{eq:two-phase-condition}
      \min\l(\sum_{i=1}^n m_i^0 \beta_i,\sum_{i=1}^n m_i^0 \frac{1}{\beta_i}\r) > 1.
  \end{equation}
  Moreover, under \eqref{eq:two-phase-condition}, this stationary state is unique and explicitly computable from $\ovX$, which is itself solution to a convex scalar equation. 
  \end{proposition}
Let us remark that, under condition \eqref{eq:two-phase-condition}, one-phase stationary states are not expected to be stable.

\section{Finite Volume Scheme}

We consider $N \in \N^*$ reference cells of uniform size $\Delta x = \frac{1}{N}$. The $N+1$ edge vertices are denoted by $0=x_{\frac{1}{2}}, x_{\frac{3}{2}},\dots,x_{N+\frac{1}{2}}=1$. We consider a time horizon $T >0$ and a time discretization with mesh parameter $\Delta t$ defined such that $N_T \Delta t = T$ with $N_T \in \N^*$. The concentrations are discretized as $\bbc_{\Delta x}^p = (c^p_{i,K})_{i \in \{1,\dots,n\},~ K \in \{1,\dots,N\}}$ for $p \in \{0,\dots,N_T\}$. The interface is discretized in time as $X^p$ for $p \in \{0,\dots,N_T\}$, and we denote by $x_{K^p+\frac{1}{2}} \in [0,1]$ the closest vertex to $X^p$ (the left vertex in case of equality). At time $t^{p-1}=(p-1)\Delta t$, the mesh is locally modified around $X^{p-1}$: the cells $K^{p-1}$ and $K^{p-1}+1$ are deformed, as presented initially in Figure~\ref{fig:displacement}, where we denote by $K$ the interface cell to alleviate the notations. To account for this deformation, we introduce $\Delta_K^{p-1}$ the size of cell $K$ at discrete time $t^{p-1}$:
\begin{equation}
\label{eq:size-cell}
\Delta_K^{p-1} = 
\begin{cases}
  (X^{p-1} - x_{K^{p-1}-\frac{1}{2}}) & \text{if } K=K^{p-1},\\
  (x_{K^{p-1}+\frac{3}{2}}-X^{p-1}) & \text{if } K=K^{p-1}+1, \\
  \Delta x & \text{otherwise.} 
\end{cases}
\end{equation}
With this notation, the initial concentrations $\bbc^0 \in L^\infty(\Omega_0; \mathbb A)$ are naturally discretized as $ c_{i,K}^0 = \frac{1}{\Delta_K^0} \int_{K} c_i^0 ~ dx $. Starting from the knowledge of $\bbc_{\Delta x}^{p-1}, X^{p-1}$, our scheme consists in
\begin{itemize}
\item[i)] solving the conservation laws and updating the interface position, leading to $(\bbc_{\Delta x}^{p,\star},X^p)$.
\item[ii)] updating the mesh to $\Delta_K^p$ and post-processing the interface concentrations into the final values $\bbc_{\Delta x}^p$.
\end{itemize} 

\subsection{Conservation Laws}
 
 The conservation laws \eqref{eq:conservation} are discretized implicitly as, for $K \in \{1,\dots,N\}, i \in \{1,\dots,n\}$, 
\begin{subequations}
\label{eq:scheme}
\begin{equation}
\label{eq:conservation-discr}
\frac{1}{\Delta t} (\Delta_{K}^{p,\star} c^{p,\star}_{i,K}- \Delta_{K}^{p-1}c_{i,K}^{p-1}) + J^{p,\star}_{i,K+\frac{1}{2}}-J^{p,\star}_{i,K-\frac{1}{2}} = 0.
\end{equation}
where we have introduced the intermediate quantity (see the intermediate mesh in Figure~\ref{fig:displacement})
\begin{equation}
  \label{eq:size-cell-semi-impl}
  \Delta_K^{p,\star} = 
  \begin{cases}
    (X^p - x_{K^{p-1}-\frac{1}{2}}) & \text{if } K=K^{p-1},\\
    (x_{K^{p-1}+\frac{3}{2}}-X^p) & \text{if } K=K^{p-1}+1, \\
    \Delta x & \text{otherwise.} 
  \end{cases}
\end{equation}
The bulk fluxes \eqref{eq:size-exclusion-flux}-\eqref{eq:Maxwell-Stefan-flux} are discretized in a way that preserves the bulk part of the dissipation structure \eqref{eq:dissip}. We refer to \cite{cances2020} (resp. \cite{cances2020d}) for the discretization of \eqref{eq:size-exclusion-flux} (resp.\eqref{eq:Maxwell-Stefan-flux}) in a single-phase and fixed domain context, since we prefer to highlight our contribution to the treatment of the interface coupling. Because of the moving interface, a correction term $-X'(t) \bbc$ appears in \eqref{eq:conservation-discr} in the interface cells, see \eqref{eq:interface-conservation}, and the numerical interface fluxes are given by a discretization of \eqref{eq:BV} as
\begin{equation}
    \label{eq:numerical-BV}
    J_{i,K^{p-1}+\frac{1}{2}}^{p,\star} = F_i^{p,\star} = c_{i,K^{p-1}}^{p,\star} \exp\left(\frac{\mu_i^{*,g}-\mu_i^{*,s}}{2} \right) - c_{i,(K^{p-1}+1)}^{p,\star} \exp \left(\frac{\mu_i^{*,s}-\mu_i^{*,g}}{2} \right),
\end{equation}
Finally, \eqref{eq:interface} is discretized as
\begin{equation}
\label{eq:interface-discrete}
X^p = X^{p-1} - \Delta t \sum_{i=1}^n F_i^{p,\star}.
\end{equation}
\end{subequations}
We denote the solution to \eqref{eq:scheme} by $(\bbc_{\Delta x}^{p,\star},X^p)$. 

\subsection{Post-Processing}
\begin{figure}
  \centering
\begin{tikzpicture}[scale=1]
\tikzmath{\xl=0; \xr = 12; \y=-4; \yf=-8;}

% 1st picture: initial 
\draw[gray, thick] (\xl,0) -- (\xr,0);
\filldraw[black] (\xl,0) circle (2pt) node[below]{0};
\filldraw[black] (\xr,0) circle (2pt) node[below]{1};
\draw[dashed] (5,0.2) -- (5,-0.2) node[below]{$x_{K+\frac{1}{2}}$}; 
\draw (3,+0.2) -- (3,-0.2) node[below]{$x_{K-\frac{1}{2}}$};
\draw (7,0.2) -- (7,-0.2) node[below]{$x_{K+\frac{3}{2}}$};
\draw (9,0.2) -- (9,-0.2) node[below]{$x_{K+\frac{5}{2}}$};
\draw[very thick,blue] (5.7,0.5) -- (5.7,-0.7) node[below]{$X^{p-1}$};
\draw [decorate,decoration={brace,amplitude=5pt,raise=4ex}]
  (3,0) -- (5.7,0) node[midway, yshift=3em]{$c_{i,K}^{p-1}$};
\draw [decorate,decoration={brace,amplitude=5pt,raise=4ex}]
  (5.7,0) -- (7,0) node[midway, yshift=3em]{$c_{i,K+1}^{p-1}$};
\draw [decorate,decoration={brace,amplitude=5pt,raise=4ex}]
  (7,0) -- (9,0) node[midway, yshift=3em]{$c_{i,K+2}^{p-1}$};
% below
\draw [decorate,decoration={brace,amplitude=5pt,raise=8ex}]
  (5.7,0) -- (3,0) node[midway, yshift=-5em]{$\Delta_K^{p-1}$};
\draw [decorate,decoration={brace,amplitude=5pt,raise=8ex}]
(7,0) -- (5.7,0) node[midway, yshift=-5em]{$\Delta_{K+1}^{p-1}$};
\draw [decorate,decoration={brace,amplitude=5pt,raise=8ex}]
  (9,0) -- (7,0) node[midway, yshift=-5em]{$\Delta_{K+2}^{p-1}$};
% Second picture: intermediate
\draw[gray, thick] (\xl,\y) -- (\xr,\y);
\filldraw[black] (\xl,\y) circle (2pt) node[below]{0};
\filldraw[black] (\xr,\y) circle (2pt) node[below]{1};
\draw[dashed] (5,\y+0.2) -- (5,\y-0.2) node[below]{$x_{K+\frac{1}{2}}$}; 
\draw (3,\y+0.2) -- (3,\y-0.2) node[below]{$x_{K-\frac{1}{2}}$};
\draw (7,\y+0.2) -- (7,\y-0.2) node[below]{$x_{K+\frac{3}{2}}$};
\draw (9,\y+0.2) -- (9,\y-0.2) node[below]{$x_{K+\frac{5}{2}}$};
%\draw[dashed] (5.7,\y+0.3) -- (5.7,\y-0.3) node[below]{$X^{p-1}$};
\draw[very thick,blue] (6.3,\y+0.5) -- (6.3,\y-0.5) node[below]{$X^{p}$};
\draw [decorate,decoration={brace,amplitude=5pt,raise=4ex}]
  (3,\y) -- (6.3,\y) node[midway, yshift=3em]{$c_{i,K}^{p,\star}$};
\draw [decorate,decoration={brace,amplitude=5pt,raise=4ex}]
  (6.3,\y) -- (7,\y) node[midway, yshift=3em]{$c_{i,K+1}^{p,\star}$};
\draw [decorate,decoration={brace,amplitude=5pt,raise=4ex}]
  (7,\y) -- (9,\y) node[midway, yshift=3em]{$c_{i,K+2}^{p,\star}$};
  % below
\draw [decorate,decoration={brace,amplitude=5pt,raise=8ex}]
(6.3,\y) -- (3,\y) node[midway, yshift=-5em]{$\Delta_K^{p,\star}$};
\draw [decorate,decoration={brace,amplitude=5pt,raise=8ex}]
(7,\y) -- (6.3,\y) node[midway, yshift=-5em]{$\Delta_{K+1}^{p,\star}$};
\draw [decorate,decoration={brace,amplitude=5pt,raise=8ex}]
  (9,\y) -- (7,\y) node[midway, yshift=-5em]{$\Delta_{K+2}^{p,\star}$};

  % Third picture: final
\draw[gray, thick] (\xl,\yf) -- (\xr,\yf);
\filldraw[black] (\xl,\yf) circle (2pt) node[below]{0};
\filldraw[black] (\xr,\yf) circle (2pt) node[below]{1};

\draw (5,\yf+0.2) -- (5,\yf-0.2) node[below]{$x_{K+\frac{1}{2}}$}; 
\draw (3,\yf+0.2) -- (3,\yf-0.2) node[below]{$x_{K-\frac{1}{2}}$};
\draw[dashed] (7,\yf+0.2) -- (7,\yf-0.2) node[below]{$x_{K+\frac{3}{2}}$};
\draw (9,\yf+0.2) -- (9,\yf-0.2) node[below]{$x_{K+\frac{5}{2}}$};
\draw[very thick,blue] (6.3,\yf+0.5) -- (6.3,\yf-0.5) node[below]{$X^{p}$};

\draw [decorate,decoration={brace,amplitude=5pt,raise=4ex}]
  (3,\yf) -- (5,\yf) node[midway, yshift=3em]{$c_{i,K}^{p}$};
\draw [decorate,decoration={brace,amplitude=5pt,raise=4ex}]
  (5,\yf) -- (6.3,\yf) node[midway, yshift=3em]{$c_{i,K+1}^{p}$};
\draw [decorate,decoration={brace,amplitude=5pt,raise=4ex}]
  (6.3,\yf) -- (9,\yf) node[midway, yshift=3em]{$c_{i,K+2}^{p}$};
  % below
  \draw [decorate,decoration={brace,amplitude=5pt,raise=8ex}]
  (5,\yf) -- (3,\yf) node[midway, yshift=-5em]{$\Delta_K^{p}$};
  \draw [decorate,decoration={brace,amplitude=5pt,raise=8ex}]
  (6.3,\yf) -- (5,\yf) node[midway, yshift=-5em]{$\Delta_{K+1}^{p}$};
  \draw [decorate,decoration={brace,amplitude=5pt,raise=8ex}]
  (9,\yf) -- (6.3,\yf) node[midway, yshift=-5em]{$\Delta_{K+2}^{p}$};
\end{tikzpicture}
\caption{A virtual mesh displacement between $t^{p-1} = (p-1) \Delta t$ and $t^p = p \Delta t$.}
\label{fig:displacement}
\end{figure}
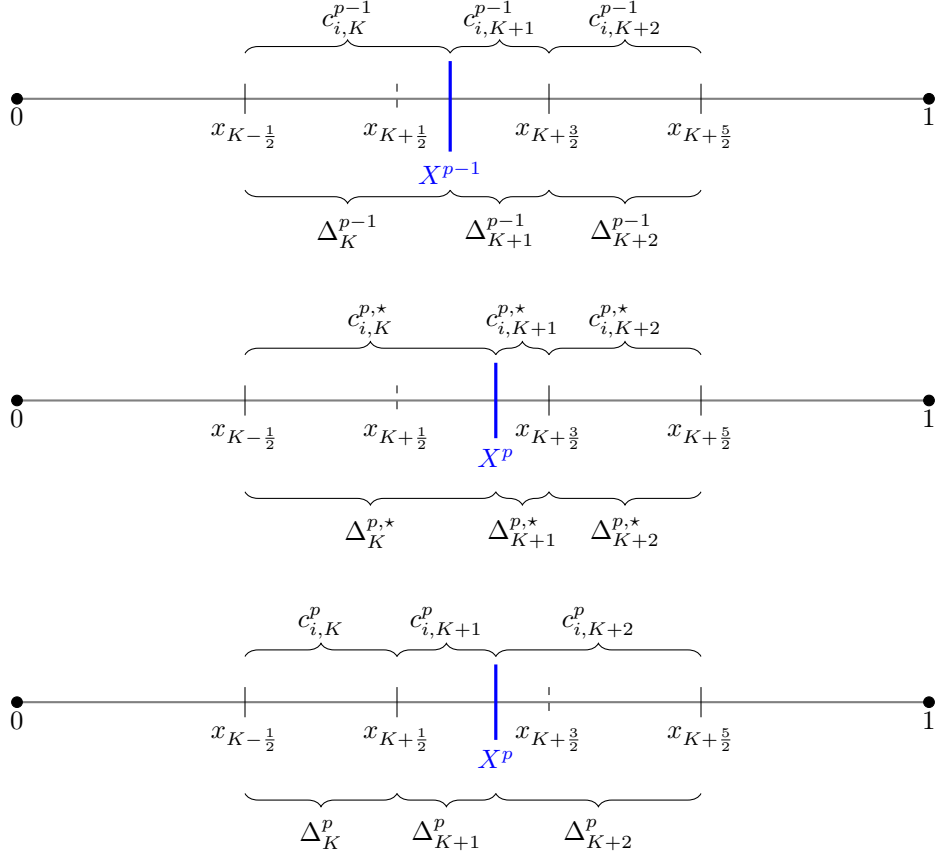

When $X^p$ crosses the center of a cell, one needs to update the interface cell from $K^{p-1}$ to $K^p$ and to adjust the concentrations accordingly. First, we can derive from \eqref{eq:interface-discrete} a linear CFL condition to enforce $|X^p-X^{p-1}| \leq \frac{\Delta x}{2}$, which in particular ensures that $|K^p-K^{p-1}| \leq 1$ and simplifies the post-processing process ($X^p$ cannot cross $x_{K+\frac{3}{2}}$ in Figure~\ref{fig:displacement}). If $K^p=K^{p-1}$, then we can directly iterate the scheme with $\bbc_{\Delta x}^p = \bbc_{\Delta x}^{p,\star}$. Otherwise, let us illustrate the case of a right displacement $K^p = K^{p-1}+1$ and let us use again the notation $K := K^{p-1}$ for simplicity. We perform the following steps (see the final mesh in Figure~\ref{fig:displacement})
\begin{itemize}
\item[i)] \emph{Projection:} The value $c_{i,K}^{p,\star}$ is assigned to the virtual cell $(x_{K-\frac{1}{2}},X^p)$. We assign this value to both the fixed cell $K=(x_{K-\frac{1}{2}},x_{K+\frac{3}{2}})$ and the new interface cell $(K+1)=(x_{K+\frac{1}{2}},X^p)$:
\begin{equation}
  \label{eq:projection}
c_{i,K}^p = c_{i,K+1}^p = c_{i,K}^{p,\star}.
\end{equation}
\item[ii)] \emph{Average:} $X^p$ replaces $x_{K+1}$ as the interface node. We average the value in the cell $(K+2)=(X^p,x_{K+2})$:
\begin{equation}
  \label{eq:average}
  c_{i,K+2}^p = \frac{1}{\Delta x + \Delta_{K+1}^{p,\star}} \left[ \Delta_{K+1}^{p,\star} c_{i,K+1}^{p,\star} + \Delta x ~ c_{i,K+2}^{p,\star} \right].
\end{equation}
\item[iii)] In all other cells, $c_{i,K}^p = c_{i,K}^{p,\star}$. 
\end{itemize} 

\subsection{Numerical Analysis}
Let us introduce the discrete version of the free energy functional \eqref{eq:free-energy}:

\begin{equation}
    \label{eq:discrete-free-energy}
    \begin{aligned}
    \cH^p(\bbc_{\Delta x}^p,X^p) &= \sum_{i=1}^n \sum_{K \leq K^p} \Delta_K^p h^s(c_{i,K}^p) + \sum_{i=1}^n \sum_{K \geq K^p+1} \Delta_K^p h^g(c_{i,K}^p).
    \end{aligned}
\end{equation}
Proposition~\ref{prop:structure} gives some a priori estimates fulfilled by a solution to the scheme, leading to existence of a solution.

\begin{proposition}[Structure preservation]
  \label{prop:structure}
Given an admissible solution \newline
$(\bbc_{\Delta x}^{p-1},X^{p-1})$, there exists an admissible solution $(\bbc_{\Delta x}^p,X^p)$ to the scheme \eqref{eq:scheme}. Moreover, the amount of matter of each species is conserved and a discrete version of the dissipation relation \eqref{eq:dissip} is satisfied:
\begin{align*}
  \bbc_{\Delta x}^p \geq 0, \; \text{and} \;  \sum_{i=1}^n c_{i,K}^p &= 1, ~ K \in \{1,\dots,N\}, \\
  \sum_{K=1}^N \Delta_K^p c_{i,K}^p &= m_i^0,~ i \in \{1,\dots,n\}, \\
  \cH^p(\bbc_{\Delta x}^p,X^p) &\leq \cH^p(\bbc_{\Delta x}^{p-1},X^{p-1}).
\end{align*}
\end{proposition}
We sketch some ingredients of the proof below, see \cite{ccce2023} for details. 

\begin{proof}
Concerning conservation of matter, it follows from summing the conservation laws \eqref{eq:conservation-discr} over the cells $K$ and the fact that the fluxes are conservative that, for any $i \in \{1,\dots,n\}$,
\[ \sum_{K=1}^N \Delta_{K}^{p,\star} c^{p,\star}_{i,K} = \sum_{K=1}^N \Delta_{K}^{p-1}c_{i,K}^{p-1}. \]
If $K^p=K^{p-1}$, the result follows immediately. Otherwise, it follows by construction of the post-processing formulas \eqref{eq:projection}-\eqref{eq:average}.  \newline
The proof of the nonnegativity of the concentrations follows from a contradiction argument with an appropriate truncation of the fluxes. One even obtains strict positivity if $\bbc_{\Delta x}^{p-1} > 0$. \newline
The volume-filling constraints are proved by summing the conservation laws \eqref{eq:conservation-discr} over $i$ and using a normalized version of \eqref{eq:interface-discrete}. \newline 
Thanks to strict positivity, a chain rule holds \cite{cances2020,cances2020d} and the continuous dissipative structure \eqref{eq:dissip} can be translated at the discrete level. Besides, convexity implies that the post-processing \eqref{eq:projection}-\eqref{eq:average} cannot make the free energy increase. \newline
Finally, the existence proof follows from a topological degree argument, arguing by deformation to two independent one-phase systems in fixed domains. 

\end{proof}

\section{Numerical Results}
The numerical scheme has been implemented in the Julia language. The nonlinear system is solved with Newton method and stopping criterion $\|Res\|_{l^2(\Delta x)} \leq 10^{-12}$, where $Res$ is the residual of the scheme. Jacobians are efficiently automatically computed thanks to the ForwardDiff and SparseDiffTools packages.  

Let us introduce a test case: we fix an initial interface $X^0=0.51$ and smooth initial concentrations $c^0_{1}(x) = c^0_{2}(x) = \frac{1}{4} \l(1 + \cos(\pi x)\r), ~  c_{3}^0(x) = \frac{1}{2}\l(1-\cos(\pi x)\r)$ that will be suitably discretized. The cross-diffusion coefficients are taken equal in each phase, with values $\kappa_{12}=\kappa_{21}= 0.2, ~ \kappa_{23}=\kappa_{32}=0.1, ~ \kappa_{13}=\kappa_{31}=1$ (diagonal coefficients do not play any role). The reference chemical potential $\bbmu^{*,s},\bbmu^{*,g}$ are given by
$e^{\bbmu^{*,s}} = [0.2 ~ 0.4 ~ 0.4], ~ e^{\bbmu^{*,g}} = [1.2 ~ 0.1 ~ 0.1]$, so as to fulfill the equilibrium condition \eqref{eq:two-phase-condition}. 

\begin{figure*}
  \centering
  \begin{subfigure}[b]{0.45\textwidth}
      \centering
      \includegraphics[width=\textwidth]{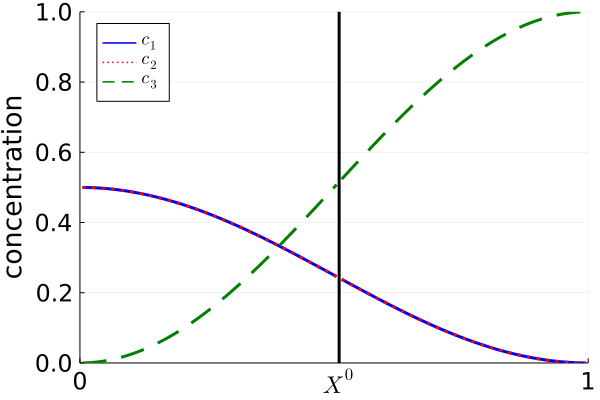}
      \caption[Network2]%
      {{\small Initial profiles}}    
  \end{subfigure}
  \hfill
  \begin{subfigure}[b]{0.45\textwidth}  
      \centering 
      \includegraphics[width=\textwidth]{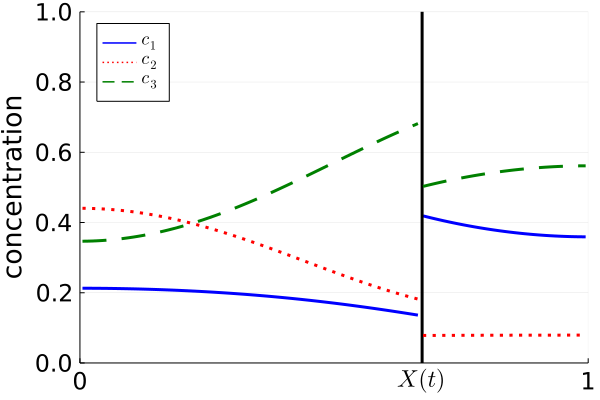}
      \caption[]%
      {{\small $t=0.25$}}    
  \end{subfigure}
  \vskip\baselineskip
  \begin{subfigure}[b]{0.45\textwidth}   
      \centering 
      \includegraphics[width=\textwidth]{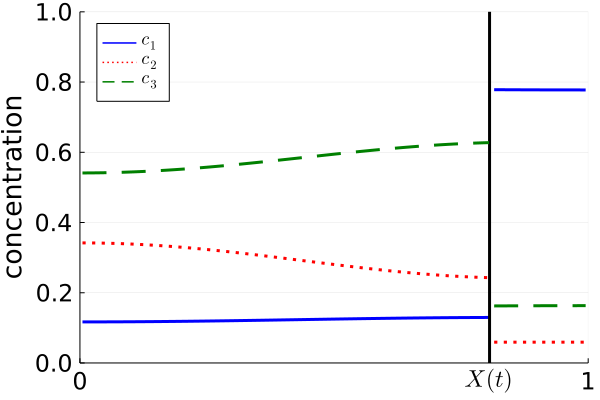}
      \caption[]
      {{\small $t=1.0$}}    
  \end{subfigure}
  \hfill
  \begin{subfigure}[b]{0.45\textwidth}   
      \centering 
      \includegraphics[width=\textwidth]{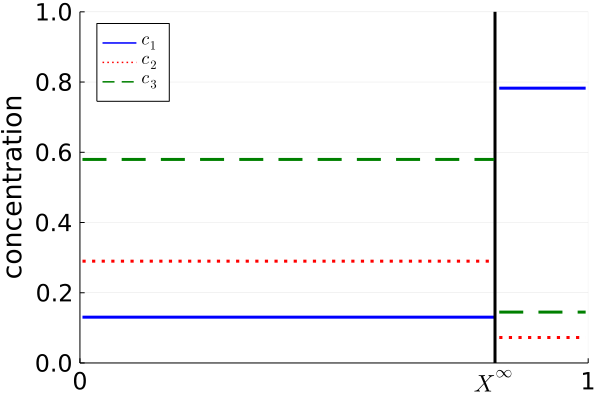}
      \caption[]%
      {{\small Stationary profiles}}  
  \end{subfigure}
  \caption[]
  {\small Concentration profiles at different times. } 
  \label{fig:snapshots}
\end{figure*}

\begin{figure*}
  \centering
  \begin{subfigure}[b]{0.45\textwidth}
    \centering
    \includegraphics[width=\textwidth]{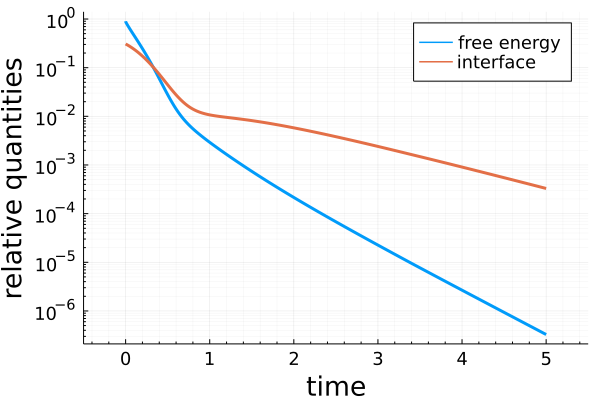}
\caption{Long-time asymptotics}
\label{fig:longtime}
\end{subfigure}
\hfill
  \begin{subfigure}[b]{0.45\textwidth}
      \centering
      \includegraphics[width=\textwidth]{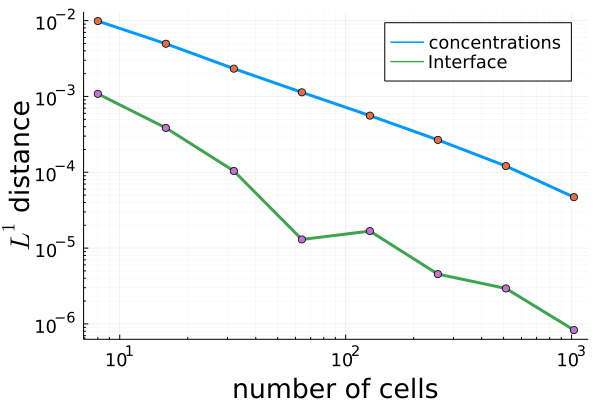}
  \caption{Convergence analysis}
  \label{fig:convergence}
  \end{subfigure} 
  \caption{\small $(\cH(\bbc(t),X(t)) - \cH(\bbc^{\infty},X^{\infty}))$ and $(X^{\infty}-X(t))$ as functions of time (left). Convergence analysis of the solution under space grid refinement (right).}
\end{figure*}

We illustrate the properties of the scheme on a uniform mesh of $N=100$ cells with time step $\Delta t_1 = 6 \times 10^{-4}$ and a final time $T_1=5$. Snapshots of the simulation are presented in Figure~\ref{fig:snapshots}, where one notices the formation of a discontinuity at the free interface and convergence to the two-phase stationary solution. We also study the long-time asymptotics: we first compute accurately the stationary solution $(\bbc^{\infty},X^{\infty})$ obtained in Proposition~\ref{prop:stationary}. Then we study the relative free energy $\cH^p(\bbc_{\Delta x}^p,X^p)-\cH^\infty(\bbc^{\infty},X^{\infty})$ and relative interface $X^{\infty}-X^p$ over time. The results are given in Figure~\ref{fig:longtime}, indicating exponential speed of convergence and decrease of both functionals. In particular, our scheme is well-balanced and preserves the asymptotics of the continuous system.

Our second test is devoted to a convergence analysis with respect to the size of the mesh. We consider a fixed time step $\Delta t_2=10^{-4}$, a final time $T_2=0.25$, uniform meshes from $2^3$ to $2^{10}$ cells and we compare the different solutions with respect to a reference solution computed on a finer grid of $2^{11}$ cells. The space-time (resp. time) $L^1$ error on the concentrations (resp. on the interface) are displayed in Figure~\ref{fig:convergence}. One clearly observes convergence, at first order in space for the concentrations. These results should be compared with the second order accurate one-phase schemes \cite{cances2020,cances2020d}. On the one hand, it is plausible that the interface treatment induces the loss of order. On the other hand, the discrete $L^1((0,1))$ space distance we use to compare solutions is not perfectly adapted since the solutions are defined in slightly different domains. Rescaling all quantities might offer more insights into the convergence properties.

\FloatBarrier

\medskip

\noindent \textbf{Acknowledgment} The authors acknowledge support from the ANR project
COMODO (ANR-19-CE46-0002) which funds the Ph.D. of Jean Cauvin-Vila.

\printbibliography

\end{document}